\documentclass [pre, showpacs] {revtex4}

\begin{document}

\title{The Mean Distance to the $n$-th Neighbour in a Uniform Distribution
of Random Points: An Application of Probability Theory}

\author{Pratip Bhattacharyya}
 \email{pratip@cmp.saha.ernet.in}

\author{Bikas K. Chakrabarti}
 \email{bikas@cmp.saha.ernet.in}

\affiliation {Theoretical Condensed Matter Physics Division,
  Saha Institute of Nuclear Physics,
  Sector - 1, Block - AF, Bidhannagar, Kolkata 700 064, India}

\date{September 17, 2003}

\begin{abstract}

We study different ways of determining the mean distance
$\left \langle r_n \right \rangle$ between a reference point
and its $n$-th neighbour among random points distributed with
uniform density in a $D$-dimensional Euclidean space.
First we present a heuristic method; though this method provides
only a crude mathematical result, it shows a simple 
way of estimating $\left \langle r_n \right \rangle$.
Next we describe two alternative means of deriving the exact expression
of $\langle r_n \rangle$: we review the method using absolute
probability and develop an alternative method using conditional
probability. Finally we obtain an approximation to
$\left \langle r_n \right \rangle$ from the mean
volume between the reference point and its $n$-th neighbour
and compare it with the heuristic and exact results.

\end{abstract}

\pacs{05.90.+m, 02.50.Cw}

\maketitle

\section {Introduction}

\indent Consider random geometrical points, i.e., points with
uncorrelated positions, distributed uniformly in a $D$-dimensional
Euclidean space, with a density of $N$ points per unit volume.
A point is said to be the $n$-th neighbour of another (the reference point)
if there are exactly $n - 1$ other points that are closer to the latter
than the former. We address the following problem:
what is the mean distance $\left\langle r_n \right\rangle$
between a given reference point and its $n$-th neighbour, $n < N$?
This is essentially a problem of geometrical interest. However
the quantity $\left\langle r_n \right\rangle$ is relevant
in certain physical and computational contexts: for example,
in astrophysics it gives the mean distance between neighbouring
stars distributed independently in a homogeneous model of the
universe~\cite{Chandrasekhar1943}. In optimization theory the
values of $\left \langle r_n \right \rangle$ help in estimating
the optimal length of a closed path connecting a given set of points
in space, as in the case of the travelling salesman
problem~\cite{Beardwood1959, Armour1983, Rujan1988, Percus1996, Cerf1997}.
It may also help in determining the statistical properties of
complex networks~\cite{Albert2002}.

\indent In this article we study different ways of determining
the mean $n$-th neighbour distance. We first present a heuristic
method of estimating $\left \langle r_n \right \rangle$ by using
a physical picture. Next we describe the derivation of the exact
expression of $\langle r_n \rangle$ in two ways:
we review the method using absolute probability and derive the
result by an alternative method using conditional probability.
The former method is comprehensive while the latter is more
analytic and explicitly illustrates the notion of ensemble average,
i.e., the mean value of a macroscopic quantity calculated over
all possible configurations of a system.
Finally we calculate the mean volume which separates the reference
point from its $n$-th neighbour and obtain an approximate expression
of $\langle r_n \rangle$ from its radius. We find out how this
approximation deviates from the exact expression of $\langle r_n \rangle$.

\section{A Heuristic Method}

\indent We begin with a heuristic approach to find the mean first
neighbour distance. Consider a unit volume of the space
described in the introduction, say,
in the form of a hypersphere or a hypercube containing exactly $N$
random points including the reference one. Let us divide this unit
volume into $N$ equal parts. Since the $N$ random points are distributed
uniformly over the unit volume each part is expected to contain just
one of these. The mean distance $\left \langle r_1 \right \rangle$
between any point and its first neighbour is naively given by the
linear extent of each part. Since the volume of each part is
$1/N$, we expect

\begin{equation}
\left \langle r_1(N) \right \rangle \approx \left ( {1 \over N} \right )
 ^{1/D}.
\label{eq:heuristic-av-1}
\end{equation}

\indent We extend the above result for the first neighbour to the
$n$-th neighbour by the following heuristic argument.
We choose any one of the points as the reference and locate its
$n$-th neighbour, $n < N$. The expected distance between them
is $\left\langle r_n(N)\right\rangle$. Keeping these two points
fixed we change the number of points in the unit volume to $N \alpha$
by adding or removing points at random; the factor $\alpha$ is arbitrary
to the extent that $N \alpha$ and $n \alpha$ are natural numbers.
Since the distribution of points is uniform, the hypersphere that
had originally enclosed $n$ points is now expected to contain
$n \alpha$ points. Therefore, what was originally the $n$-th neighbour
of the reference point is now expected to be the $n \alpha$-th neighbour
for which the expected distance from the reference point is now
$\langle r_{n \alpha} \left ( N \alpha \right ) \rangle$.
Since the two points under consideration are fixed, so is the distance
between them. Consequently

\begin{equation}
\left\langle r_n(N)\right\rangle
 \approx \left\langle r_{n \alpha}(N \alpha)\right\rangle .
\label{eq:heuristic-step1}
\end{equation}

\noindent The above relation is approximate as the change in the density
of points does not always convert the $n$-th neighbour of the
reference point to exactly its $n \alpha$-th neighbour.
Now we take $\alpha = 1/n$, so that

\begin{equation}
\left\langle r_n(N)\right\rangle
 \approx \left\langle r_1 \left ( {N \over n} \right ) \right\rangle
\label{eq:heuristic-step2}
\end{equation}

\noindent which shows that the mean $n$-th neighbour distance for
a set of $N$ random points distributed uniformly is approximately
given by the mean first neighbour distance for a depleted set
of $N/n$ random points in the same volume. The above relation is
derived for such values of $n$ that divide $N$ exactly; however this
approximate relation may be used for any value of $n$ when $N \gg n$,
by replacing $N/n$ with the integer nearest to it. Using the expression
of $\left\langle r_1(N)\right\rangle$ from Eq.~(\ref{eq:heuristic-av-1})
we get

\begin{equation}
\left \langle r_n(N) \right \rangle \approx \left ( {n \over N} \right )
 ^{1/D}
\label{eq:heuristic-av-n}
\end{equation}

\noindent which, therefore, requires as a necessary condition that
$N \gg n$. The results of Eq.~(\ref{eq:heuristic-av-1}) and
Eq.~(\ref{eq:heuristic-av-n}) are extremely crude approximations;
however these provide us with a rough picture of the mean $n$-th
neighbour distance.

\section{The Method of Absolute Probablity}

\indent We now review the derivation of the exact expression of
$\left \langle r_n \right \rangle$ by using the theory of
{\it absolute probability}~\cite{Feller1968}. The method described
is similar to that followed in Ref.~\cite{Cerf1997}.
Consider the system of random points described at the begining
of the introduction. Assuming a certain random point as the
reference there will be $N - 1$ other random points within a
$D$-dimensional hypersphere of unit volume with the reference point at
its center. For a given reference point the absolute probability
of finding its $n$-th neighbour ($n < N$)
at a distance between $r_n$ and $r_n + {\rm d}r_n$ from it is
given by the probability that out
of the $N - 1$ random points (other than the reference point)
distributed uniformly within the hypersphere of unit volume,
exactly $n - 1$ points lie within a concentric hypersphere of radius
$r_n$ and at least one of the remaining $N - n$ points lie within
the shell of internal radius $r_n$ and thickness ${\rm d}r_n$:

\begin{equation}
P(r_n) \: {\rm d}r_n =
     \left ( \begin{array}{c}
                      N - 1\\
                      n - 1
             \end{array} \right ) \: V_n^{n - 1} \:
     \sum_{q=1}^{N-n} \left ( \begin{array}{c}
                                                 N - n\\
                                                 q
                                        \end{array} \right )
               \left ( 1 - V_n \right )^{N - n - q} \:
               \left ( {\rm d}V_n \right ) ^q
\label{eq:abs-prob-n-orig}
\end{equation}

\noindent where
\begin{equation}
V_n = {\pi^{D/2} \over \Gamma \left ( {D \over 2} + 1 \right ) } \: r_n^D 
\label{eq:D-dim-vol}
\end{equation}

\noindent is the volume of the $D$-dimensional hypersphere of radius
$r_n$ centered at the reference point and

\begin{equation}
\left ( \begin{array}{c}
                 N - 1\\
                 n - 1
        \end{array} \right ) = {(N - 1)! \over (n - 1)! \: (N - n)!},
\hspace{1.0cm}
\left ( \begin{array}{c}
                 N - n\\
                 q
        \end{array} \right ) = {(N - n)! \over q! \: (N - n - q)!}
 \nonumber
\label{eq:binomial-coeff}
\end{equation}

\noindent are binomial coefficients. Ignoring differentials of order
 higher than the first ($q = 1$) in Eq.~(\ref{eq:abs-prob-n-orig}) we get:

\begin{equation}
P(r_n) \: {\rm d}r_n = \left ( \begin{array}{c}
                                           N - 1\\
                                           n - 1
                               \end{array} \right ) \: (N - n) \:
 V_n^{n - 1} \left ( 1 - V_n \right )^{N - n - 1} \: {\rm d}V_n.
\label{eq:abs-prob-n}
\end{equation}

\noindent Since the $n$-th neighbour $(n < N)$ must certainly
lie within a unit volume centered at the reference point, its mean
distance from the reference point is given by:

\begin{equation}
\left \langle r_n \right \rangle = \int_0^R \: r_n \: P(r_n) \:
 {\rm d}r_n,
\label{eq:abs-defn-av-n}
\end{equation}

\noindent where $R$ is the radius of the $D$-dimensional hypersphere
of unit volume:

\begin{equation}
R = {\left [\Gamma \left ({D \over 2} + 1 \right ) \right ]^{1/D}
 \over \pi^{1/2}}.
\label{eq:unitvol-rad}
\end{equation}

\noindent Changing the variable of integration in
Eq.~(\ref{eq:abs-defn-av-n}) from radius to volume (by the relation of
Eq.~(\ref{eq:D-dim-vol})) and using the probability distribution of
Eq.~(\ref{eq:abs-prob-n}) we get the exact result for the mean $n$-th
neighbour distance:

\begin{eqnarray}
\left\langle r_n(N)\right\rangle & = &
 {\left [ \Gamma \left ( {D \over 2} + 1 \right ) \right ]^{1/D}
  \over \pi^{1/2}} \:
                            \left ( \begin{array}{c}
                                                   N - 1\\
                                                   n - 1
                                    \end{array} \right ) (N - n) \:
  \int_0^1 V_n^{n + (1/D) - 1} \:
 \left ( 1 - V_n \right ) ^{N - n - 1} \: {\rm d}V_n \nonumber \\
 & = & { \left [ \Gamma \left ( {D \over 2} + 1 \right ) \right ]^{1/D}
 \over \pi^{1/2}} \:
                                  \left ( \begin{array}{c}
                                                    N - 1\\
                                                    n - 1
                                          \end{array} \right )
 \: (N - n) \: B \left ( n + {1 \over D}, N - n \right ) \nonumber \\
 & = & { \left [ \Gamma \left ( {D \over 2} + 1 \right ) \right ]^{1/D}
 \over \pi^{1/2}} \:
 {\Gamma \left ( n + {1 \over D} \right ) \over \Gamma(n)}
  \: {\Gamma(N) \over \Gamma \left ( N + {1 \over D} \right ) }.
\label{eq:exact-av-n}
\end{eqnarray}

\noindent Here $B(x, y)$ is the beta function defined as
$B(x, y) = \int_0^1 t^{x - 1} (1 - t)^{y - 1} {\rm d}t$
and $\Gamma(z)$ is the complete gamma function:
$\Gamma(z) = \int_0^\infty t^{z-1} e^{-t} {\rm d}t$.
These functions are related by the formula~\cite{Abramowitz1972}~:
$B(x, y) = \Gamma(x) \: \Gamma(y) / \Gamma(x + y)$.

\indent For large values of $N$ we get by using Stirling's
approximation~\cite{Graham1994} for the gamma function:
$\Gamma \left ( N + 1/D \right ) /
\Gamma \left ( N \right ) \sim N^{1/D}$;
therefore, for a large density $N$, Eq.~(\ref{eq:exact-av-n})
reduces to the following asymptotic form:

\begin{equation}
\langle r_n(N) \rangle \sim
 { \left [ \Gamma\left( {D \over 2} + 1 \right ) \right ]^{1/D}
 \over \pi^{1/2}} \:
 \: {\Gamma \left ( n + {1 \over D} \right ) \over \Gamma(n)}
     \left ( {1 \over N} \right )^{1/D}.
\label{eq:asympN-av-n}
\end{equation}

\noindent If the neighbour index $n$ is also large (but
$n < N$), we have $\Gamma \left ( n + 1/D \right ) /
\Gamma \left ( n \right ) \sim n^{1/D}$ and the complete
asymptotic expression of the mean $n$-th neighbour distance
is given by:

\begin{equation}
\langle r_n(N) \rangle \sim
 { \left [ \Gamma\left( {D \over 2} + 1\right ) \right ]^{1/D}
 \over \pi^{1/2}} \:
 \left ( {n \over N} \right )^{1/D}.
\label{eq:asymp-av-n}
\end{equation}

\noindent The above equation shows that the expression of
$\left \langle r_n(N) \right \rangle$ obtained by heuristic means
[Eq.~(\ref{eq:heuristic-av-n})] has the correct asymptotic dependence
on $N$ and $n$.

\section{The Method of Conditional Probability}

\indent Next we develop an alternative way of deriving the
exact expression of $\left \langle r_n \right \rangle$ by
using the theory of {\em conditional probability}~\cite{Feller1968}.
We proceed by asserting that we look for the $n$-th neighbour
of a reference point only after its first $n - 1$ neighbours have
been located. In that case the reference point and its first $n - 1$
neighbours are considered as given; now the probablity
${\cal P}(r_n) {\rm d}r_n$ of finding the $n$-th neighbour of
the reference point at a distance between $r_n$ and
$r_n + {\rm d}r_n$ from it is a {\em conditional probability}
as the $n$-th neighbour must certainly lie outside the hypersphere
of radius $r_{n-1}$:

\begin{equation}
{\cal P}\left (r_n\right ) \: {\rm d}r_n =
 \sum_{q=1}^{N-n} \left ( \begin{array}{c}
                                 N - n\\
                                 q
                          \end{array} \right ) \:
 \left [ 1 - {V_n - V_{n - 1}
 \over 1- V_{n - 1}}\right ]^{N - n - q} \:
 \left [ {{\rm d}V_n \over 1 - V_{n - 1}} \right ]^q .
 \label{eq:cond-prob-n-orig}
\end{equation}

\noindent Here $V_n$ is the volume of the $D$-dimensional hypersphere
of radius $r_n$ centered at the reference point [Eq.~(\ref{eq:D-dim-vol})].
Ignoring differentials of order higher than the first ($q = 1$)
in Eq.~(\ref{eq:cond-prob-n-orig}) we get:

\begin{equation}
{\cal P}\left (r_n\right ) \: {\rm d}r_n =
 \left [ 1 - {V_n - V_{n - 1} \over 1- V_{n - 1}}\right ]^{N - n - 1}
 \: {(N - n) \: {\rm d}V_n \over 1 - V_{n - 1}} .
 \label{eq:cond-prob-n}
\end{equation}

\indent For a given reference point and its first $n - 1$ neighbours
the mean $n$-th neighbour distance is thus obtained as:

\begin{equation}
\left \langle r_n \right\rangle^{(\rm conditional)}
 = \int_{r_{n-1}}^R \: r_n \: {\cal P}\left ( r_n\right) \: {\rm d}r_n
 \label{eq:cond-defn-av-n}
\end{equation}

\noindent where, as before, $R$ is the radius of a $D$-dimensional
hypersphere of unit volume. The quantity 
$\left \langle r_n\right\rangle^{(\rm conditional)}$
is a function of the particular $r_{n-1}$, $r_{n-2}$, $\ldots$, $r_1$
which are the distances of the first $n-1$ neighbours of the
reference point. To remove the dependence of the mean $n$-th
neighbour distance on the particular set of values of the first
$n - 1$ neighbour distances the quantity
$\left \langle r_n\right\rangle^{(\rm conditional)}$
must be averaged successively over the probability distributions
of each of the first $n-1$ neighbours:

\begin{equation}
\left \langle r_n \right\rangle =
   \int_0^R \: {\rm d}r_1 \: {\cal P}(r_1) \:
   \int_{r_1}^R \: {\rm d}r_2 \: {\cal P}(r_2) \:
   \cdots \: \int_{r_{n-3}}^R \: {\rm d}r_{n-2} \: {\cal P}(r_{n-2})
   \: \int_{r_{n-2}}^R \: {\rm d}r_{n-1} \: {\cal P}(r_{n-1}) \:
   \left \langle r_n\right\rangle^{(\rm conditional)}
 \label{eq:interm-av-n}
\end{equation}

\noindent where the probability distribution of the $i$-th neighbour
distance is given by Eq.~(\ref{eq:cond-prob-n}) with $i$ replacing $n$.
This step is equivalent to an {\em ensemble average} in statistical
mechanics~\cite{Kittel1958}. After changing all the variables of
integration in Eq.~(\ref{eq:interm-av-n})
from radii to the corresponding volumes (by  the relation of
Eq.~(\ref{eq:D-dim-vol})) and changing the order of the integrals such that
the integral with respect to $V_n$ has to be evaluated last, we get:

\begin{eqnarray}
\left \langle r_n(N) \right\rangle & = &
 {\left [\Gamma\left ( {D \over 2} + 1\right )\right ]^{1/D}
 \over \pi^{1/2}} \: (N - 1) (N - 2) \cdots (N - n) \int_0^1 \: {\rm d}V_n \:
 V_n^{1/D} \: \left (1 - V_n \right )^{N - n - 1} \nonumber \\
 & & \times \: \int_0^{V_n} \: {\rm d}V_1 \:
               \int_{V_1}^{V_n} \: {\rm d}V_2 \: \cdots \:
     \int_{V_{n-3}}^{V_n} \: {\rm d}V_{n-2} \:
              \int_{V_{n-2}}^{V_n} \: {\rm d}V_{n-1}
 \label{eq:changed-general}
\end{eqnarray}

\noindent which gives the final form of the mean $n$-th neighbour
distance:

\begin{eqnarray}
\left \langle r_n(N) \right \rangle & = &
{\left [ \Gamma \left ( {D \over 2} + 1\right ) \right ]^{1/D}
 \over \pi^{1/2}}
      \left ( \begin{array}{c}
                       N - 1\\
                       n - 1
              \end{array} \right ) \: (N - n) \: \int_0^1 V_n^{n + (1/D) - 1}
     \left (1 - V_n \right )^{N - n - 1} \: {\rm d}V_n \nonumber \\
 & = & { \left [ \Gamma \left ( {D \over 2} + 1 \right ) \right ]^{1/D}
         \over \pi^{1/2}} \:
       {\Gamma \left ( n + {1 \over D} \right ) \over \Gamma(n)} \:
       {\Gamma(N) \over \Gamma \left ( N + {1 \over D} \right ) }.
 \label{eq:exact-av-n'}
\end{eqnarray}

\noindent This is identical to the result of Eq.~(\ref{eq:exact-av-n})
obtained in the previous section.

\section {The mean volume estimate}

\indent Instead of calculating the mean distance to the $n$-th neighbour
we now calculate the mean volume $\left \langle V_n \right \rangle$
separating the reference point from its $n$-th neighbour.
The volume separating a reference point from its $n$-th neighbour
located at a distance $r_n$ from it is defined as the volume of the
hypersphere of radius $r_n$ and centered at the reference point.
Therefore, from Eq.~(\ref{eq:D-dim-vol}), we get:

\begin{equation}
\left \langle V_n \right \rangle = {\pi^{D/2} \over
 \Gamma \left ( {D \over 2} + 1 \right ) }
 \left \langle r_n^D \right \rangle .
\label{eq:mean-vol-n-defn}
\end{equation}

\noindent Using the absolute probability distribution of
Eq.~(\ref{eq:abs-prob-n}) we get

\begin{equation}
\left \langle r_n^D \right \rangle
 = \int_0^R r_n^D \: P(r_n) \: {\rm d}r_n
\label{eq:mean-r-Dmoment}
\end{equation}

\noindent where $R$ is defined in Eq.~(\ref{eq:unitvol-rad}). Consequently,
we get from the above two equations:

\begin{eqnarray}
\left \langle V_n \right \rangle
   & = & \left ( \begin{array}{c}
                           N - 1 \nonumber \\
                           n - 1
                 \end{array} \right ) \:
         (N - n) \int_0^1 V_n^n \:
         \left ( 1 - V_n \right )^{N - n - 1} \: {\rm d}V_n
         \nonumber \\
 ~ & = & \left ( \begin{array}{c}
                           N - 1 \nonumber \\
                           n - 1
                 \end{array} \right ) \: (N - n) \: B(n+1, N-n) \nonumber \\
 ~ & = & {\Gamma(n+1) \over \Gamma(n)} {\Gamma(N) \over \Gamma(N+1)}
         = {n \over N}.
\label{eq:mean-vol-n}
\end{eqnarray}

\indent Now we may estimate the mean $n$-th neighbour distance by
using the following approximation:

\begin{equation}
\left \langle r_n \right \rangle
 \approx \left \langle r_n^D \right \rangle^{1/D}.
\label{eq:mean-vol-est}
\end{equation}

\noindent From Eqs.~(\ref{eq:mean-vol-n-defn}) and (\ref{eq:mean-vol-n})
we get:

\begin{equation}
\langle r_n^D(N) \rangle^{1/D} =
 { \left [ \Gamma \left( {D \over 2} + 1 \right ) \right ]^{1/D}
 \over \pi^{1/2}}
 \left ( {n \over N} \right )^{1/D}.
\label{eq:approx-av-n}
\end{equation}

\noindent The above result is the same as the asymptotic expression of
$\left \langle r_n \right \rangle$ obtained in Eq.~(\ref{eq:asymp-av-n}) 
and thus~Eq.(\ref{eq:mean-vol-est}) is valid as $N \to \infty$
and $n \to \infty$, $n < N$.
Comparison with Eq.~(\ref{eq:heuristic-av-n}) shows that this approximation
is a better estimate of $\langle r_n \rangle$ than the result obtained
by heuristic means.

\indent The error in this estimate of $\left \langle r_n \right \rangle$
by Eq.~(\ref{eq:mean-vol-est}) is given by
$\left \langle r_n^D \right \rangle^{1/D} - \left \langle r_n \right \rangle$.
From Eqs.~(\ref{eq:exact-av-n}) and (\ref{eq:approx-av-n}) we get:

\begin{eqnarray}
{\pi^{1/2} \over
 \left [ \Gamma \left( {D \over 2} + 1 \right ) \right ]^{1/D}}
 \left [ \left \langle r_n^D \right \rangle^{1/D}
 - \left \langle r_n \right \rangle \right ]
 & = & \left ( {n \over N} \right )^{1/D}
       - \: {\Gamma \left ( n + {1 \over D} \right ) \over \Gamma(n)} \:
      {\Gamma(N) \over \Gamma \left ( N + {1 \over D} \right ) }.
\label{eq:approx-av-n-error}
\end{eqnarray}

\noindent It is obvious that the error is zero for $D = 1$ and
from the expression in Eq.~(\ref{eq:approx-av-n-error})
it is clear that the error is greater than zero for all finite
dimensions $D \geq 2$. For large values of $D$ we get:

\begin{equation}
{\Gamma \left ( n + {1 \over D} \right ) \over \Gamma \left ( n \right )}
 \approx \Gamma \left ( 1 + {1 \over D} \right ) \:
 \left [ 1 + {H_{n-1} \over D} \right ]
\label{eq:gamma-ratio-large-D}
\end{equation}

\noindent and a similar expression for $\Gamma \left ( N + 1/D \right ) /
\Gamma \left (N \right )$, where $H_{n} = \sum_{k=1}^n 1/k$ are called
harmonic numbers~\cite{Graham1994}. Consequently
Eq.~(\ref{eq:approx-av-n-error}) reduces to the form:

\begin{eqnarray}
{\pi^{1/2} \over
 \left [ \Gamma\left( {D \over 2} + 1 \right ) \right ]^{1/D}}
 \left [ \left \langle r_n^D \right \rangle^{1/D}
 - \left \langle r_n \right \rangle \right ]
 & \approx & \left ( {n \over N} \right )^{1/D} -
 {1 + {1 \over D} H_{n-1} \over 1 + {1 \over D} H_{N-1}} \nonumber \\
 & \approx & \left ( {n \over N} \right )^{1/D}
   - \left [ 1 + {1 \over D} \left ( H_{n-1} - H_{N-1} \right ) \right ].
\label{eq:error-large-D}
\end{eqnarray}

\noindent Since $\lim_{n \to \infty} (H_n - \ln n) = \gamma =
\lim_{N \to \infty} (H_N - \ln N)$, where $\gamma = 0.5772156649 \ldots$
is the Euler's constant~\cite{Graham1994}, for large values of $n$ and $N$
Eq.~(\ref{eq:error-large-D}) may be written as:

\begin{equation}
{\pi^{1/2} \over
 \left [ \Gamma\left( {D \over 2} + 1 \right ) \right ]^{1/D}}
 \left [ \left \langle r_n^D \right \rangle^{1/D}
 - \left \langle r_n \right \rangle \right ]
 \approx \left ( {n \over N} \right )^{1/D}
   - \left [ 1 + {1 \over D} \ln \left ( {n \over N} \right ) \right ].
\label{eq:error-large-DNn}
\end{equation}

\noindent The error is then practically zero, since, for large values
of $D$, we get:

\begin{equation}
\left ( {n \over N} \right )^{1/D}
 = \exp{\left [ {1 \over D} \ln \left ( {n \over N} \right ) \right ]}
 \approx 1 + {1 \over D} \ln \left ( {n \over N} \right ).
\label{eq:why-error-zero}
\end{equation}

\section{Concluding remarks}

\indent In this article we have studied two kinds of approaches
to determine the mean $n$-th neighbour distance in a system of
uniformly distributed random points. In one kind of approach we
construct the solution from a physical picture of the system;
though it produces only an approximate mathematical
result it helps to visualise the solution. The heuristic method of
section 2 and the mean volume method of section 5 are of this kind.
The other kind of approach is rigorous and produces the exact
mathematical result~: while the method of absolute probability
used in section 3 is largely pedagogical, the method of conditional
probability used in section 4 provides a detailed insight into
the problem.
Though the heuristic estimate of $\langle r_n \rangle$ is close
to the exact result only for large values of $n$, $N$ and $D$,
along with the condition $N \gg n$, the advantage of the heuristic
method lies in its simplicity; however approximate,
it gives an essence of $\langle r_n \rangle$. Therefore, in cases
where the distribution of points in space is such that an exact
evaluation of $\langle r_n \rangle$ is not possible~\cite{Percus1998},
heuristic constructions similar to this one may be useful.

\begin{acknowledgments}

We thank D. Dhar, A. G. Percus, S. S. Manna and A. Chakraborti
for their comments.

\end{acknowledgments}

\end{document}